\documentclass[12pt]{amsart}
\usepackage{amssymb}
\usepackage[dvips]{graphics}
\ifx\mathrm\undefined\let\mathrm\rm\fi
\ifx\mathbf\undefined\let\mathbf\bf\fi
\ifx\mathfrak\undefined\let\mathfrak\frak\fi
\ifx\mathcal\undefined\let\mathcal\cal\fi
\ifx\mathbb\undefined\let\mathbb\Bbb\fi
\ifx\emph\undefined\let\emph\it\fi
 at9.98pt

\newcommand{\g}{{{\mathfrak g}\,}}

\newcommand{\n}{{{\mathfrak n}}}
\newcommand{\B}{{{\mathfrak b}}}
\newcommand{\h}{{{\mathfrak h\,}}}

\newcommand{\Z}{{\mathbb Z}}

\newcommand{\C}{{\mathbb C}}

\newcommand{\Ref}[1]{{(\ref{#1})}}

\newcommand{\la}{\lambda}

\newcommand{\dontprint}[1]
{\relax}

\newtheorem%
{thm}{Theorem}[section]
\newtheorem%
{proposition}[thm]{Proposition}
\newtheorem%
{lemma}[thm]{Lemma}
\newtheorem%
{lemmadef}[thm]{Lemma-Definition}
\newtheorem%
{corollary}[thm]{Corollary}
\newtheorem%
{conjecture}[thm]{Conjecture}

\newcommand{\bea}{\begin{eqnarray*}}
\newcommand{\eea}{\end{eqnarray*}}
\newcommand{\bean}{\begin{eqnarray}}
\newcommand{\eean}{\end{eqnarray}}

\newcommand{\nc}{\newcommand}
\nc{\on}{\operatorname}
\nc{\al}{\alpha}
\nc{\ri}{\rangle}
\nc{\lef}{\langle}
\nc{\W}{{\mathcal W}}
\nc{\La}{\Lambda}
\nc{\ep}{\epsilon}
\nc{\Om}{\Omega}
\newcommand{\be}{\begin{displaymath}}
\newcommand{\ee}{\end{displaymath}}
\newcommand{\bs}{\boldsymbol}
\nc{\PCr}{{ \Bbb P  (\C[x])^r   }}
\newtheorem{theorem}{Theorem}[section]
\newtheorem{conj}[theorem]{Conjecture}

\newcommand{\p}{{\partial_{h}}}

\newcommand{\bP}{{\Bbb P^1}}

\newcommand{\id}{{\rm{id}}}
\newcommand{\om}{{\rm{Om}_{D_{\bs y}}}}

\newcommand{\tN}{{\rm{Om}_{D_{\bs y^0}}}}

\begin{document}

\title[Discrete Miura Opers and Solutions of the Bethe Ansatz
Equations]
{Discrete Miura Opers and Solutions of the Bethe Ansatz Equations}

\author[{}]{ Evgeny Mukhin  
${}^{*,1}$ 
\and Alexander Varchenko${}^{**,2}$}
\thanks{${}^1$ Supported in part by NSF grant DMS-0140460}

\thanks{${}^2$ Supported in part by NSF grant DMS-0244579}

\begin{abstract}
Solutions of the Bethe ansatz equations associated to the XXX model of
a simple Lie algebra $\g$
come in families called the populations. We prove
that a population is isomorphic to the flag variety of 
the Langlands dual Lie algebra ${}^t\g$.
The proof is based on the correspondence between the solutions of the Bethe 
ansatz equations and special difference operators which we call 
the discrete Miura opers. The notion of a discrete Miura oper is 
one of the main results of the paper.

For  a discrete Miura oper $D$,  associated to a point of a
population, we show that all solutions of the difference equation $DY=0$
are rational functions, and the solutions 
can be written explicitly in terms of points composing the population.

\end{abstract}

\maketitle 
\medskip \centerline{\it ${}^*$
Department of Mathematical Sciences,}
\centerline{\it Indiana University
Purdue University Indianapolis,}
\centerline{\it 402 North Blackford St., Indianapolis,
IN 46202-3216, USA}
% \newline mukhin@math.iupui.edu}
 \medskip
\centerline{\it ${}^{**}$Department of Mathematics, University of
  North Carolina at Chapel Hill,} \centerline{\it Chapel Hill, NC
  27599-3250, USA} \medskip

%\maketitle

\centerline{January, 2004}

\section{Introduction}

The Bethe ansatz is a large collection of methods in the theory of
quantum integrable models to calculate the spectrum and eigenvectors
for a certain commutative subalgebra of observables for an integrable
model. This commutative subalgebra includes the Hamiltonian of the
model. Its elements are called integrals of motion or conservation
laws of the model. Most of recent developments of the Bethe ansatz
methods is due to the quantum inverse scattering transform, invented
by the Leningrad school of mathematical physics. The bibliography on
the Bethe ansatz method is enormous. We refer the reader to reviews
\cite{BIK, Fa, FT}.

In the theory of the Bethe ansatz one assigns the Bethe ansatz equations to
an integrable model. Then a solution of the Bethe ansatz equations gives an
eigenvector of commuting Hamiltonians of the model.  The general
conjecture is that the constructed vectors form a basis in the space of 
states of the model. The first step to that conjecture is
to count the number of solutions of the Bethe ansatz equations.
One can expect that the number of solutions is equal to the dimension
of the space of states of the model.

The Bethe ansatz equations of the XXX model is a system of algebraic equations
associated to a Kac-Moody algebra $\g$, a non-zero step $h \in \C$,
complex numbers $z_1,\dots,z_n$, integral dominant
$\g$-weights $\La_1,\dots,\La_n$ and an integral $\g$-weight $\La_\infty$, 
see \cite{OW} and Section \ref{Bethe eqn sec with param} below.  

To approach the corresponding counting problem, to
every solution of the Bethe ansatz equations
we assign  an object called the
population of solutions. We expect that it would be easier to count populations than
individual solutions. For instance, if the Kac-Moody algebra is of type $A_r$,
then each population corresponds to a point of the intersection
of suitable Schubert varieties in a suitable Grassmannian variety, as was
shown in \cite{MV3}. 
Then the Schubert calculus allows us to count the number of intersection points
of the Schubert varieties and to give an upper bound on the number 
of populations, see \cite{MV1} and \cite{MV3}.

A population of solutions is an interesting object. It is an
algebraic variety. It is finite-dimensional, if the Weyl group of
the Kac-Moody algebra is finite. In this
paper we prove that a population is isomorphic to the flag variety of
the Langlands dual Lie algebra $^t\g$.  The proof is based on the
correspondence between solutions of the Bethe ansatz equations
 and special difference operators which we call the discrete Miura opers.

Let $G$ be the complex simply connected group with Lie algebra $^t\g$.
To every solution $\bs t$ of the Bethe ansatz equations we assign
a linear difference operator $D_{\bs t}= \p - V_{\bs t}$
where $\p : f(x) \to f(x+h)$ is the shift operator and
$V_{\bs t}(x)$ is a suitable rational $G$-valued function.
We call that difference operator a
 discrete Miura oper. Our discrete Miura opers are
analogs of the differential operators considered by V. Drinfeld and
V. Sokolov in their study of the KdV type equations \cite{DS}.

Different solutions of the Bethe ansatz equations
correspond to different discrete Miura opers. 
The discrete Miura opers, corresponding to points of a given population,
form an equivalence class with respect to a suitable gauge
equivalence. Thus a population is isomorphic to an equivalence class
of discrete Miura opers. 
We show that an equivalence class of discrete Miura opers is
isomorphic to the flag variety of $^t\g$.

\medskip

If $D_{\bs t}$ is the discrete Miura oper corresponding to a solution $\bs t$
of the Bethe ansatz equations, then the set of solutions of the difference equation
$D_{\bs t}\,Y\,=\,0$ with values in a suitable space is an important
characteristics of $\bs t$.  It turns out that, for any simple Lie algebra and
any solution $\bs t$ of the Bethe ansatz equations, 
the difference equation $D_{\bs t}\,Y\,=\,0$  has a rational fundamental matrix 
of solutions.  Moreover, all
solutions of the difference equation $D_{\bs t}\,Y\,=\,0$ can be
written explicitly in terms of points composing the
population, originated at $\bs t$.  Thus,
 the population of solutions of the Bethe ansatz equations
``solves''  the Miura difference equation $D_{\bs t}\,Y\,
=\,0$ in rational functions. This is the second main result of the paper.

\medskip

The results of this paper in the cases of $A_r$ and $B_r$ were first obtained
in \cite{MV3}. 

\medskip

The populations related to the Gaudin model of a Kac-Moody algebra
 $\g$ were introduced in \cite{MV1}. In \cite{MV1} we conjectured that
 every $\g$-population is isomorphic to the flag variety of the
 Langlands dual Kac-Moody algebra $\g^L$.  That conjecture
was proved  for $\g$ of type  $A_r, B_r, C_r$ in \cite{MV2},
for $\g$ of type $G_2$ in \cite{BM},
for all simple Lie algebras in \cite{F1} and \cite{MV4}. 
The ideas of \cite{MV4} motivated the present paper.

\medskip

There are different versions of the XXX Bethe ansatz equations
associated to a simple Lie algebra, see \cite{OW, MV2, MV3}.
Ogievetsky and Wiegman introduced in \cite{OW} a set of Bethe ansatz
equations for any simple Lie algebra $\g$. For $\g$ of type $A_r, D_r,
E_6, E_7, E_8$ the Ogievetsky-Wiegman equations are the Bethe ansatz
equations considered in this paper. For other
simple Lie algebras the Ogievetsky-Wiegman equations are different
from the Bethe ansatz equations considered in this paper, see Section
\ref{Sym par}.

\medskip

The discrete Miura opers, considered in this paper, are 
discrete versions of the special differential operators
called the Miura opers and introduced in \cite{DS}. 
The Miura opers play an essential role
in the Drinfeld-Sokolov reduction and geometric Langlands correspondence,
see \cite{DS, FFR, F2, FRS}. It would be interesting to see if our discrete opers
may play a similar role in discrete versions of the Drinfeld-Sokolov reduction
and geometric Langlands correspondence.

Considerations of the present paper
are in the spirit of the geometric Langlands correspondence.
Namely, we start from a solution $\bs t$
of the Bethe ansatz equations associated to a simple Lie algebra $\g$, that is
we start from a Bethe
 eigenvector of the commuting Hamiltonians of the XXX
$\g$ model. Having a solution $\bs t$ we construct
the associated discrete Miura $G$-oper $D_{\bs t}$, whose fundamental matrix
is a rational function. Having 
$D_{\bs t}$ we may recover the $\g$ population of solutions of the Bethe ansatz
equations originated at $\bs t$. In particular we may recover $\bs t$ and the
associated Bethe eigenvector of the commuting Hamiltonians.

The fact that the discrete Miura oper has a rational fundamental matrix of
solutions is a discrete analog of the fact that a differential Miura oper
has the trivial monodromy group, see \cite{F1, MV4}. 

\medskip

The paper is organized as follows. In Section 2 we introduce populations 
of solutions of the Bethe ansatz equations. 
In Section 3 we discuss elementary properties of discrete Miura opers
corresponding to solutions of the Bethe ansatz equations. 
In Section 4 we give explicit formulas for 
solutions of the difference equation $D_{\bs t}\,Y\,=\,0$, see Theorems
 \ref{sl},  \ref{Solutions}.
In Section 5 we prove that the variety of 
gauge equivalent marked discrete
Miura $\g$ opers is isomorphic to the flag variety of $^t\g$, see Theorem
\ref{isomorphism}. We discuss the relations between the Bruhat cell decomposition
of the flag variety and the populations of solutions of the Bethe ansatz equations
in Section 6. The main results of the paper are Corollaries 
\ref{ISO} and \ref{main corollary}.

\medskip

We thank P. Belkale and S. Kumar for numerous useful discussions.

\section{Populations of Solutions 
of the Bethe Equations, \cite{MV2}}\label{crit pts}

\subsection{Kac-Moody algebras}\label{Kac_Moody sec}
Let $A=(a_{i,j})_{i,j=1}^r$ be a generalized  Cartan matrix, 
$a_{i,i}=2$,
$a_{i,j}=0$ if and only $a_{j,i}=0$,
 $a_{i,j}\in \Z_{\leq 0}$ if $i\ne j$. 
We  assume that $A$ is symmetrizable, i.e. 
there exists a diagonal matrix $D=\on{diag}\{d_1,\dots,d_r\}$ 
with positive integers $d_i$ such that $B=DA$
is symmetric.

Let $\g=\g(A)$ be the corresponding complex Kac-Moody 
Lie algebra (see \cite{K}, \S 1.2), 
$\h \subset \g$  the Cartan subalgebra.
The associated scalar product is non-degenerate on $\h^*$ and 
 $\on{dim}\h = r + 2d$,  where $d$ is the dimension of 
the kernel of the Cartan matrix $A$.

Let $\al_i\in \h^*$, $\al_i^\vee\in \h$, $i = 1, \dots , r$, be the sets of simple roots,
coroots, respectively. We have
\bea
(\al_i,\al_j)&=& d_i \ a_{i,j}, \\
\langle\la ,\al^\vee_i\rangle&=&2(\la,\al_i)/{(\al_i,\al_i)}, \qquad \la\in\h^*.
\eea
In particular, $\langle\al_j ,\al^\vee_i \rangle = a_{i,j}$.

Let $\mathcal P = \{ \lambda \in \h^* \, |\, \langle\la
,\al^\vee_i\rangle \in \Z\}$ and $\mathcal P^+ = \{ \lambda \in \h^* \, |\, \langle\la
,\al^\vee_i\rangle \in \Z_{\geq 0}\}$ be the sets of integral and
dominant integral weights.

Fix $\rho\in\h^*$ such that $\langle\rho,\al_i^\vee\rangle=1$,
$i=1,\dots,r$. We have $(\rho,\al_i)= (\al_i,\al_i)/2$.

The Weyl group $\mathcal W\in\on{End (\h^*)}$ is generated by 
reflections $s_i$, $i=1,\dots,r$, 
\be
s_i(\la)=\la-\langle\la,\al_i^\vee\rangle\al_i, \qquad \la\in\h^*.
\ee
We use the notation
\bea\label{shifted}
w\cdot\la=w(\la+\rho)-\rho,\qquad w\in \mathcal W,\;\la\in\h^*,
\eea
for the shifted action of the Weyl group.

The Kac-Moody algebra $\g(A)$ is generated by $\h$, $e_1, \dots , e_r, f_1,
\dots , f_r$ with defining relations
\bea
[e_i, f_j] & = & \delta_{i,j} \,\alpha_i^\vee ,
\qquad i, j = 1, \dots r ,
\\
{}[ h , h'] & = & 0 ,
\qquad h, h' \in \h ,
\\
{}[ h, e_i] &=& \langle \alpha_i, h \rangle\, e_i ,
\qquad  h \in \h, \ i = 1, \dots r ,
\\
{}[ h, f_i] &=& - \langle \alpha_i, h \rangle\, f_i ,
\qquad  h\in \h, \ i = 1, \dots r ,
\eea
and the Serre's relations
\bea
(\mathrm{ad}\,{} e_i)^{1-a_{i,j}}\,e_j = 0 ,
\qquad
(\mathrm{ad}\, {} f_i)^{1-a_{i,j}}\,f_j = 0 ,
\eea
for all $i\neq j$. The generators
$\h$, $e_1, \dots , e_r, f_1,
\dots , f_r$ are called the Chevalley generators.

Denote  $\n_+$ (resp. $\n_-$) the subalgebra
generated by $e_1, \dots , e_r$ (resp. $f_1, \dots , f_r$). Then
$\g = \n_+\oplus \h\oplus \n_-$. Set $\B_\pm = \h \oplus \n_\pm$.

The Kac-Moody algebra $^t\g=\g( ^tA)$ 
corresponding to the transposed Cartan matrix $^tA$ is called
{\it Langlands dual} to $\g$. 
Let $^t\al_i\in {}^t\h^*$, $^t\al_i^\vee\in {}^t\h$, 
$i = 1, \dots , r$, be the sets of simple roots,
coroots of $^t\g$, respectively. Then 
\bea
\langle {}^t\al_i, {}^t\al_j^\vee\rangle  =  \langle\al_j,\al_i^\vee\rangle = a_{i,j}
\eea
 for all $i,j$.

\subsection{The Bethe ansatz equations with parameters $\bs b$}
\label{Bethe eqn sec with param}
Fix a Kac-Moody algebra
$\g=\g(A)$, a non-negative integer $n$, a collection of dominant integral weights  
$\bs\La = (\La_1, \dots , \La_n)$, $\La_i\in\mathcal P^+$,
and complex numbers $\bs z = (z_1, \dots , z_n)$.
Fix a non-zero complex number $h$.

Fix a collection  $\bs b = (b_{i,m})_{i,m = 1,\ i\neq m}^r$ of complex numbers.
We say that the parameters $\bs b$ are {\it symmetric} if they satisfy the 
condition:
\bean\label{b condition}
b_{i,m}\ +\ b_{m,i}\ =\ h
\eean
for all $i, m \in \{1, \dots , r\},\ i\neq m$.

Choose 
a collection of non-negative integers  $\bs l=(l_1,\dots,l_r)\in\Z^r_{\geq 0}$. 
The choice of $\bs l$ is equivalent to the choice of the weight 
\be
\La_\infty\ =\sum_{i=1}^n \La_i  -  \sum_{j=1}^r
l_j\al_j \in \mathcal P .
\ee
The weight $\La_\infty$ will be called {\it the weight at infinity}.
Set $\bar {\bs \Lambda} = (\La_1, \dots , \La_n, \La_\infty)$.
Let 
\bea
\bs  t\ =\  \{\ t_j^{(i)} \in \C \ | \ i = 1,\dots,r, \ j = 1, \dots , l_i\ \}\ 
\eea 
be a collection of complex numbers.

{\it The {\rm XXX} Bethe ansatz equations 
associated to} $ \bar{\bs \La},\ \bs z, \ \bs b$ is the
following system of algebraic equations with respect to the variables $\bs t$:
\bean\label{Bethe with parameters}
&&
\prod_{s=1}^n\ \frac{t_j^{(i)}\ -\ z_s\ +\ (\La_s,\al_i^\vee)\, h/2}
{t_j^{(i)}\ -\ z_s\ -\ (\La_s,\al_i^\vee)\, h/2}\ {} \times
\phantom{aaaaaaaaaaaaaaaaaaaaaaaaaaaaaaaaaa}
\\ 
&&
\ {}\
\prod_{m=1,\ m\neq i}^{r} \left(\prod_{k=1}^{l_m}
\ \frac{t_j^{(i)}\ -\ t_k^{(m)}\ +\ b_{i,m}}
{t_j^{(i)}\ -\ t_k^{(m)}\ +\ b_{i,m}\ -\  h}\right)^{-a_{i,m}}\ 
\prod_{k=1,\ k\neq j}^{l_i}
\frac{t_j^{(i)}\ -\ t_k^{(i)}\ -\ h}
{t_j^{(i)}\ -\ t_k^{(i)}\ +\ h}\
=\ 1\ ,
\notag
\eean
where $i = 1, \dots , r$,\ $j = 1, \dots , l_i$.

The product of symmetric groups
$S_{\bs l}=S_{l_1}\times \dots \times S_{l_r}$ acts on the  set of
solutions of \Ref{Bethe with parameters}
permuting the coordinates with the same upper index. 

For $i = 1, \dots , r$, consider the $l_i$ equations \Ref{Bethe with parameters}
with fixed upper index $i$. We call that system of equations
  {\it the Bethe ansatz equations with fixed  upper index $i$}.

\subsection{ Polynomials representing solutions of the 
Bethe ansatz equations }\label{PLCP}
For a given  $\bs t = (t^{(i)}_j)$ introduce an $r$-tuple of polynomials 
$\bs y=( y_1(x),$ $\dots ,$ $ y_r(x))$, where
\bean\label{y}
y_i(x)\ =\ \prod_{j=1}^{l_i}(x-t_j^{(i)}).
\eean
Each polynomial is considered up to multiplication 
by a non-zero number.
The tuple defines a point in the direct product 
$\PCr$ of $r$ copies of the projective space associated to the vector 
space of polynomials in $x$. 
We say that {\it the tuple $\bs y$ represents the collection of numbers } $\bs t$.

It is convenient to think that if  a polynomial $y_k$ of a tuple
$\bs y = (y_1, \dots , y_r) \in \PCr$ has degree zero, then it means that 
the collection $\bs t = (t^{(i)}_j)$ has no $t^{(i)}_j$-s with $i=k$.

%the tuple $(1, \dots , 1)$ of constant polynomials
%represents in $\PCr$ the collection $\bs t = (t^{(i)}_j)$
%with no $t^{(i)}_j$-s. 
%This corresponds to the case when $\bs l = (0, \dots , 0)$ and $\La_\infty =
%\sum_{s=1}^n \La_s$.

For $i = 1, \dots , r$ introduce polynomials 
\bean\label{T}
T_i(x)\ =\ \prod_{s=1}^n\prod_{p=0}^{(\La_s, \al_i^\vee) - 1}
( x\ -\ z_s\ +\ (\La_s, \al_i^\vee )\ h/2\ - \ p h )\ 
\eean
and
\bean\label{Q}
Q_i(x)\ =\ T_i(x)\ \prod_{m=1,\ m\neq i}^r \ y_m(x + b_{i,m})^{-a_{i,m}}\ .
\eean

We say that the tuple $\bs y $ is {\it generic with respect 
to weights  $\bs \La$,   numbers
$\bs z$, and parameters} $\bs b$ if 
for every $i = 1, \dots , r$  the polynomial $y_i(x)$ has no multiple roots
and no common roots with polynomials
$\ y_i(x \ + \ h),$ and  $Q_i(x)$.

If $\bs y$ represents a solution $\bs t$ of the Bethe ansatz equations 
\Ref{Bethe with parameters} and $\bs y$ is generic, then the $S_{\bs l}$-orbit of 
$\bs t$ is called {\it a Bethe solution of 
\Ref{Bethe with parameters}}.

The Bethe ansatz equations can be written as
\bean\label{Bethe with param}
\frac{Q_i(\,t^{(i)}_j\,)}{Q_i(\,t^{(i)}_j\ -\ h\,)}\
\prod_{k=1,\ k\neq j}^{l_i}
\frac{t_j^{(i)}\ -\ t_k^{(i)}\ -\ h}
{t_j^{(i)}\ -\ t_k^{(i)}\ +\ h}\
=\ 1\ ,
\eean
where $i = 1, \dots , r$,\ $j = 1, \dots , l_i$.

\bigskip

For $i = 1, \dots , r$, a tuple $\bs y$  is called {\it fertile in the $i$-th
direction with respect to $\bs \La, \ \bs z,\ \bs b$}, if there 
exists a polynomial $\tilde y_i$ satisfying the equation
\bean\label{wron}
 y_i(x + h)\ \tilde y_i(x)\  -\   y_i(x)\
\tilde y_i(x + h)\ =\ Q_i(x)\ .
\eean
A tuple $\bs y$ is called {\it fertile with respect to $\bs \La, \ \bs z,\ \bs b$}, 
if it is fertile in all directions $i = 1, \dots , r$.

\medskip

{\bf Example.} The tuple $(1, \dots , 1)$ is fertile with respect to any
given $\bs \La, \ \bs z,\ \bs b$.

\medskip

Instead of saying that $\bs y$ is generic or fertile with respect to
$\bs \La, \ \bs z,\ \bs b$ we will also say that $\bs y$ is 
{\it generic or fertile
with respect to polynomials $T_1(x), \dots , T_r(x)$ and parameters $\bs b$}.

If $\bs y$ is fertile in the $i$-th direction and $\tilde y_i$ is a polynomial
solution of \Ref{wron}, then the tuple
\bean\label{simple}
\bs y^{(i)}\ =\ (y_1, \dots , \tilde y_i,\dots, y_r) \ 
\ \in \ \  \PCr 
\eean
is called {\it an immediate descendant} of $\bs y$
in the $i$-th direction.

If $\tilde y_i$ is a solution of \Ref{wron}, then $\tilde y_i + c y_i$ is a solution too 
for any $c \in \C$.

\begin{lemma}[\cite{MV1, MV2}]
Assume that $\bs y$ is generic. Let $\tilde y_i$ be a solution of \Ref{wron}. Then the tuple
$(y_1, \dots , \tilde y_i + c y_i, \dots , y_r)$ is generic for almost all $c\in \C$. The 
exceptions form a finite subset in $\C$.
\hfill
$\square$
\end{lemma}

\begin{lemma}[\cite{MV1, MV2}]
Let $\bs y^j$, $j = 1, 2, \dots $, be a sequence of tuples in $\PCr$ which has a limit
$\bs y$. Assume that all tuples $\bs y^j$ are fertile. Then $\bs y$ is fertile.
\hfill
$\square$
\end{lemma}

\begin{lemma}[see \cite{MV1, MV2}] \label{reflction}
Denote $\tilde l_i = \deg \tilde y_i$
and $\La_\infty^{(i)} = \sum_{s=1}^n \La_s - \tilde l_i \al_i -
\sum_{j=1, \ {} j\neq i}^r l_j \al_j$.
If $\tilde l_i \neq l_i$,  then 
$$
\La_\infty^{(i)}\  =\
s_i\cdot \La_\infty \ ,
$$
 where $s_i\cdot$ is the shifted action of the 
$i$-th reflection of the Weyl group.
\hfill
$\square$
\end{lemma}

\begin{thm}[cf. \cite{MV2}]\label{fertile cor}
${}$

\begin{enumerate}
\item[(i)]
Let a tuple $\bs y = (y_1, \dots , y_r)$
be generic. Let $i \in \{1, \dots , r\}$.
Then $\bs y$ is fertile in the $i$-th direction if and only if $\bs t$ satisfies
the Bethe ansatz equations \Ref{Bethe with parameters} with fixed upper index $i$.
\item[(ii)] Let parameters $\bs b$ be symmetric. 
Let $\bs y$ be generic and fertile. Let $i \in \{1, \dots , r\}$.
Let $\bs y^{(i)}$ be an immediate descendant of $\bs y$ in the $i$-th direction.
Assume that $\bs y^{(i)}$ is generic. Then $\bs y^{(i)}$ is fertile.
\end{enumerate}
\end{thm}

\begin{proof}
To prove (i) introduce $g(x) = \tilde y_i(x)/y_i(x)$ and write \Ref{wron}
as
\bean\label{g}
g(x) \ -\ g(x + h)\ = \frac{Q_i(x)}{ y_i(x)\ y_i(x+h)}\ .
\eean
The tuple $\bs y$ is fertile in the $i$-th direction if and only if there exists a rational
function $g(x)$ satisfying \Ref{g}. A function $g(x)$ exists if and only if
\bean\label{res}
\text{Res}_{x= t_j^{(i)}}\
\frac{Q_i(x)}{ y_i(x)\ y_i(x+h)}\ =\ -\
\text{Res}_{x= t_j^{(i)}-\,h}\
\frac{Q_i(x)}{ y_i(x)\ y_i(x+h)}\ 
\eean
for $j = 1, \dots , l_j$.
The systems of equations \Ref{res} is equivalent to the system of the Bethe ansatz equations
with fixed upper index $i$. Part (i) is proved.

%Let  $\bs y^{(i)}$ be an immediate descendant of $\bs y$, deg $\tilde y_i = \tilde l_i$.
%Let $\tilde t^{(i)}_j$, $j = 1, \dots , \tilde l_i$ be roots of $\tilde y_i$.

To prove (ii) we check that the Bethe ansatz equations \Ref{Bethe with parameters} 
are satisfied for roots of polynomials composing the tuple $\bs y^{(i)}$.

The Bethe ansatz equations with upper index $i$ are satisfied for roots of $\bs y^{(i)}$
according to part (i).

If $m$ is such that $a_{i,m} = 0$, then the Bethe ansatz equations with upper index $m$ for 
roots of
$\bs y^{(i)}$ are the same as for $\bs y$, since the roots of $\tilde y_i$ do not enter 
those equations.

Let $m$ be such that $a_{i,m} \neq 0$. Write \Ref{wron} as
\bea
\frac{y_i(x+h)}{y_i(x)}\ -\ \frac{\tilde y_i(x+h)}{\tilde y_i(x)}\
= \ \frac{Q_i(x)}{\tilde y_i(x)\ y_i(x)}\ .
\eea
Substitute  to this equation
the  zeros of the polynomial $y_m(x+b_{i,m})$ and get
\bean\label{zeros}
\frac{y_i(t_k^{(m)}-b_{i,m}+h)}{y_i(t_k^{(m)}-b_{i,m})}\ 
= \ \frac{\tilde y_i(t_k^{(m)}-b_{i,m}+h)}{\tilde y_i(t_k^{(m)}-b_{i,m})}\
\eean
for $k = 1, \dots , l_m$.

The Bethe ansatz equations with upper index $m$
for roots of $\bs y^{(i)}$  contain the factor
$ \tilde y_i(t^{(m)}_k + b_{m,i})/\tilde y_i(t^{(m)}_k + b_{m,i}-h)$ while
the Bethe ansatz equations with upper index $m$ for roots of $\bs y$
contain the factor
$ y_i(t^{(m)}_k + b_{m,i})/ y_i(t^{(m)}_k + b_{m,i}-h)$.
By \Ref {zeros} the two ratios are equal if parameters $\bs b$ are symmetric.
Hence the Bethe ansatz equations with upper index $m$ are satisfied for roots of
$\bs y^{(i)}$ if they are satisfied for roots of $\bs y$. Part (ii) is proved.
\end{proof}

\subsection{Simple reproduction procedure}    
Assume that parameters $\bs b$ are symmetric.

Let $\bs y$ represent a Bethe solution of \Ref{Bethe with parameters}.
Let $i \in \{1, \dots , r\}$, and let $\tilde y_i$ be a polynomial solution
of equation \Ref{wron}. For complex numbers $c_1$ and $c_2$, not both equal to zero,
consider the tuple
\bea
\bs y_{(c_1:c_2)}^{(i)}\ =\ (y_1, \dots ,\  c_1\tilde y_i + c_2 y ,\
 \dots,\ y_r)\ {}\ \in \PCr\ .
\eea 
The tuples form a one-parameter family. The parameter space of the family
is the projective line $\Bbb P^1$ with projective coordinates 
$(c_1 : c_2)$.
We have a map
\bean\label{map}
Y_{\bs y, i}\ : \ \Bbb P^1 \ \to \PCr\ , 
\qquad 
(c_1 : c_2)\ {} \mapsto \ {}
\bs y_{(c_1:c_2)}^{(i)}\ .
\notag
\eean
Almost all tuples $\bs y_{(c_1:c_2)} ^{(i)}$ are generic.
The exceptions form a finite set in $\Bbb P^1$.

Thus, starting with a tuple $\bs y$,
 representing a Bethe solution of equations \Ref{Bethe with parameters}
associated to numbers $z_1, \dots , z_n$,
integral dominant weights $\La_1, \dots , \La_n$,  a weight $\La_\infty$
at infinity, parameters $\bs b$,
 and an index $i \in \{1, \dots , r\}$, we  construct  a  family 
$Y_{\bs y, i} : \Bbb P^1 \to \PCr$ of fertile tuples.
For almost all $c \in \Bbb P^1$ (with finitely many exceptions only), 
the tuple $Y_{\bs y, i}(c)$  represents 
a Bethe solution of the Bethe ansatz equations 
associated to points $z_1, \dots , z_n$,
integral dominant weights $\La_1, \dots , \La_n$, parameters $\bs b$, 
 and a suitable weight at infinity.

We call this construction  {\it the
simple reproduction procedure in the $i$-th direction}.

\subsection{General reproduction procedure}\label{general procedure}
Assume that parameters $\bs b$ are symmetric.

Assume that a tuple $\bs y \in \PCr$ represents a 
Bethe solution of the Bethe ansatz equations 
associated to $\bs z,\ \bar{\bs \La},\ \bs b$.

Let $\bs i = (i_1,  \dots , i_k), \ 1 \leq i_j \leq r,$ be a sequence of natural numbers.
We define a $k$-parameter family of fertile tuples 
\bea\label{general map}
Y_{\bs y, \bs i} \ :\ (\Bbb P^1)^k \ \to \PCr
\eea
by induction on $k$,  starting at $\bs y$ and successively applying 
the simple reproduction procedure in directions
$i_1, \dots , i_k$. The image of this map is denoted by
$P_{\bs y, \bs i}$ .

For a given $\bs i = (i_1, \dots , i_k)$, almost all tuples 
$Y_{\bs y, \bs i} ( \bs c)$ represent  
Bethe solutions of the Bethe ansatz equations associated to 
 points $z_1, \dots , z_n$, dominant integral weights
 $\La_1, \dots , \La_n$, symmetric parameters $\bs b$, and suitable weights at infinity.
Exceptional values of $c \in (\Bbb P^1)^k$ are contained in a proper algebraic subset.

It is easy to see that if $\bs i' = (i'_1,  \dots , i'_{k'}),
 \ 1 \leq i'_j \leq r$,
is a sequence of natural numbers, 
and the sequence $\bs i'$ is contained in the sequence $\bs i$ as  an ordered subset, then
$P_{\bs y, \bs i'}$ is a subset of $P_{\bs y, \bs i}$.

The union 
\be 
P_{\bs y} \ = \ \cup_{\bs i} \ P_{\bs y, \bs i}\ \subset \PCr\ ,
\ee
where the summation is over all of sequences $\bs i$, 
is called {\it the population of solutions of the Bethe ansatz equations associated} 
to the Kac-Moody algebra $\g$, integral dominant
weights $\La_1, \dots , \La_n$, numbers $z_1, \dots , z_n$, 
symmetric parameters $\bs b$, and {\it originated} at $y$.

If two populations with the same $\bs \La, \ \bs z, \ \bs b$
intersect, then they coincide.

If the Weyl group is finite, then all tuples of a population consist of polynomials 
of bounded degree. Thus, if the Weyl group of $\g$ is finite, then a population
is an irreducible projective variety.

Every population $P$ has a tuple $\bs y = (y_1, \dots , y_r)$,\ $\deg\,y_i = l_i$, such that
the weight $\La_\infty = \sum_{s=1}^n \La_s - \sum_{i=1}^r l_i \al_i$ is dominant integral,
see \cite{MV1, MV2}.

\begin{conj}[\cite{MV2}]\label{CON}
Every population, associated to a Kac-Moody algebra $\g$, dominant integral weights
$\Lambda_1, \dots , \Lambda_n$,   points $z_1, \dots , z_n$,
symmetric parameters $\bs b$,
is an algebraic variety isomorphic to 
the flag variety associated to the Kac-Moody algebra
$^t\g$ which is Langlands dual to $\g$. Moreover, the parts of the family corresponding to 
tuples of polynomials with fixed degrees are isomorphic to 
Bruhat cells of the flag variety.
\end{conj}

The conjecture is proved  for the Lie algebra of type $A_r$ in \cite{MV3}. 
In this paper  we prove the conjecture for every simple Lie algebra.

%For $A_r$ the proof of the conjecture
%in this paper is different from the proof in \cite{MV2}.

\subsection{Special symmetric parameters $\bs b$}\label{Sym par}
Here are two examples of symmetric parameters $\bs b$.

The parameters $\bs b$ are symmetric if
\bean\label{sym 0}
b_{i,m} \ = \frac h 2 \
\qquad
\text{for\ all} \ i\neq m\ .
\eean
The parameters $\bs b$ satisfying \Ref{sym 0} will be called {\it the 
Ogievetsky-Wiegman parameters}, cf. \cite{OW, MV3}.

The parameters $\bs b$ are symmetric if
\bean\label{sym}
b_{i,m} \ = \ 0 \ {}\ {}
\text{for} \ i > m
\qquad
\text{and} 
\qquad
b_{i,m} \ = \ h \ {}\ {}
\text{for} \ i < m\ .
\eean
The parameters $\bs b$ satisfying \Ref{sym} will be called {\it the special
symmetric parameters}, cf. \cite{OW, MV2}.

Let $\bs b^1 = (b_{i,m}^1)_{i,m = 1,\ i\neq m}^r$ and
$\bs b^2 = (b_{i,m}^2)_{i,m = 1,\ i\neq m}^r$ be two collections of parameters.
We say that $\bs b^1$ and $\bs b^2$ are {\it gauge equivalent} if there exist
complex numbers $d^{(1)}, \dots , d^{(r)}$ with the following property.
We require that for any 
tuple $(y_1(x), \dots , y_r(x))$, fertile respect to some polynomials $T_1(x), \dots ,
T_r(x)$ and parameters $\bs b^1$, the tuple
$(y_1(x+d^{(1)}), \dots , y_r(x+d^{(r)}))$ is fertile with 
respect to polynomials $T_1(x + d^{(1)}), \dots , T_r(x + d^{(r)})$ and
parameters $\bs b^2$.

If $\bs b^1$ and $\bs b^2$ are gauge equivalent and
$P$ is a population  associated to
some polynomials $T_1(x), \dots , T_r(x)$ and parameters $\bs b^1$,
then the set 
\bea
\{ (y_1(x+d^{(1)}), \dots , y_r(x+d^{(r)})) \ | \ (y_1(x), \dots ,
y_r(x)) \in P \}
\eea
is a population associated to polynomials $T_1(x+d^{(1)}), \dots , T_r(x+d^{(r)})$
and parameters $\bs b^2$.

\begin{thm}\label{b equiv}
Let $\bs b^1$ and $\bs b^2$ be symmetric parameters. 
Assume that the Dynkin diagram  of the Cartan matrix $A$ of the Lie
algebra $\g$ is a tree. Then $\bs b^1$ and $\bs b^2$ are gauge equivalent.
\end{thm}

\begin{proof} 
We will introduce parameters $\bs b^3 = (b^3_{i,j})$ 
in terms of the Dynkin diagram and will
prove that $\bs b^1$ and $\bs b^3$ are equivalent. That will prove the theorem.

Let $v_1, \dots , v_r$ be vertices
of the Dynkin diagram corresponding to
the roots $\al_1, \dots , \al_r$, respectively. 
For $i = 2, \dots , r$, let  
$v_{i_1}, \dots , v_{i_k}$ be the unique sequence of distinct vertices 
of the Dynkin diagram such that
for $j=1, \dots , k-1$ the vertices  $v_{i_j}$ and $v_{i_{j+1}}$ are connected
by an edge, and $v_{i_1} = v_1$, $v_{i_k} = v_i$. The number $k$ will be 
called the distance between $v_1$ and $v_i$ and denoted by $\delta_i$.
Let $\bs b^3$ be defined by the rule:
\bea
b^3_{i,j} = 0\ {} \
\text{if}\ {} \  \delta_i > \delta_j ,
\qquad
b^3_{i,j} = h \ {} \
\text{if}\ {} \  \delta_i < \delta_j ,
\eea
\bea
b^3_{i,j} = 0\ {} \
\text{if}\ {} \  \delta_i = \delta_j \ {}\
\text{and}\ {}\ i > j ,
\qquad
b^3_{i,j} = h \ {} \
\text{if}\ {} \  \delta_i = \delta_j \ {}\
\text{and}\ {}\ i < j .
\eea
Clearly $\bs b^3$ is symmetric. 

Define $(d^{(1)}, \dots , d^{(r)})$. Set
\bea
d^{(i)}\ =\ b^1_{i_1,i_2} \ +\ b^1_{i_2,i_3}\ +\ \dots\ +\
b^1_{i_{k-1},i_k}\ -\ (k-1)h \ 
\eea
for $i>1$ and set $d^{(1)} = 0$.
It is easy to see that the sequence $(d^{(1)}, \dots , d^{(r)})$ establishes
the equivalence of $\bs b^1$ and $\bs b^3$. 
\end{proof}

\begin{corollary}\label {sym par 0}
If the Dynkin diagram is a tree, then any set of symmetric parameters is  gauge
equivalent to the set of special symmetric parameters given by \Ref{sym}.
\end{corollary}

\begin{corollary}\label {sym par}
If $\g$ is simple, then any set of symmetric parameters is  gauge
equivalent to the set of special symmetric parameters given by \Ref{sym}.
\end{corollary}

Ogievetsky and Wiegman considered in \cite{OW} a set of Bethe ansatz equations for any simple
Lie agebra $\g$. For $\g$ of type $A_r, D_r, E_6, E_7, E_8$ the  Ogievetsky-Wiegman equations are
the Bethe ansatz equations associated to parameters given by \Ref{sym 0}. For other simple
Lie algebras the Ogievetsky-Wiegman equations are different from the Bethe ansatz equations 
considered in this paper.

For $\g$ of type $A_r$ 
we considered in \cite{MV3} the Bethe ansatz equations associated to
the special symmetric parameters given by \Ref{sym}.

\subsection{Diagonal sequences of polynomials associated to a Bethe solution
and a sequence of indices} In this section 
we assume that parameters $\bs b$ are symmetric.
We introduce notions which will be used
in Chapter \ref{solutions} to construct solutions of difference equations.

\begin{lemma}\label{generation}
Assume that a tuple $\bs y \in \PCr$ represents a 
Bethe solution of the Bethe ansatz equations associated to $\bs z,\ \bar{\bs \La}$,
symmetric parameters $\bs b$. 
Let $\bs i = (i_1, \dots , i_k), \ 1 \leq i_j \leq r,$ 
be a sequence of natural numbers. Then there exist tuples
$\bs y^{(i_1)} = ( y_1^{(i_1)}, \linebreak
\dots , y_r^{(i_1)})$,
$\bs y^{(i_1, i_2)} = ( y_1^{(i_1,i_2)}, \dots , y_r^{(i_1,i_2)})$,
\dots ,
$\bs y^{(i_1, \dots , i_k)} = ( y_1^{(i_1, \dots , i_k)}, \dots , 
y_r^{(i_1, \dots , i_k)})$ in \linebreak
$\PCr$ such that
\begin{enumerate}
\item[(i)] 
\bea
&&
y_{i_1}^{(i_1)}(x)\
y_{i_1}(x + h)\ -\
 y^{(i_1)}_{i_1}(x + h)\ y_{i_1}(x)\
\\
&&
\phantom{aaaa aaaa aaaa aaaa aaaa }\ 
 =\  T_{i_1}(x)\ \prod_{m=1,\ m\neq i_1}^{r}
 \  ( y_m(\, x + b_{{i_1},m}\,) )^{-a_{i_1, m}} 
\eea
and $ y_j^{(i_1)} = y_j$ for $j\neq i_1$;
\item[(ii)] 
for $l = 2, \dots , k$, 
\bea
&&
y_{i_l}^{(i_1, \dots , i_{l})}(x)\
y_{i_l}^{(i_1, \dots , i_{l-1})}(x + h)\ - \
y_{i_l}^{(i_1, \dots , i_{l})}(x + h)\ 
y_{i_l}^{(i_1, \dots , i_{l-1})}(x)
\\
&&
\phantom{aaaa aaaa aaaa aaaa}\ 
 =\  T_{i_l}(x)\ \prod_{m=1,\ m\neq i_l}^r
 \  ( y_m^{(i_1, \dots , i_{l-1})}(\, x + b_{{i_l},m}\,) )
^{- \, a_{{i_l},m}} 
\eea
and $ y_{j}^{(i_1, \dots , i_{l})} =  y_{j}^{(i_1, \dots , i_{l-1})}$
for $j\neq i_l$.
\hfill
$\square$
\end{enumerate}
\end{lemma}

The tuples 
$\bs y^{(i_1)}$, $\bs y^{(i_1, i_2)}$, \dots ,
$\bs y^{(i_1,  \dots , i_k)}$
belong to the population  $P_{\bs y}$.
The tuple $\bs y^{(i_1)}$ is obtained from $\bs y$ by the $i_1$-th
simple generation procedure and  
for $l = 2, \dots , k$ the tuple $\bs y^{(i_1, \dots , i_{l})}$
is obtained from $\bs y^{(i_1, \dots , i_{l-1})}$
by the $i_l$-th simple generation procedure.

The sequence of tuples $\bs y^{(i_1)}$, $\bs y^{(i_1,i_2)}$, \dots ,
$\bs y^{(i_1,  \dots , i_k)}$ satisfying Lemma \ref
{generation} will be called {\it associated to the Bethe solution $\bs y$
and the sequence of indices $\bs i$}. The sequence of polynomials
$y_{i_1}^{(i_1)}$, $y_{i_2}^{(i_1,i_2)}$, \dots ,
$y_{i_k}^{(i_1,  \dots , i_k)}$
 will be called {\it the diagonal sequence of polynomials associated to
the Bethe solution $\bs y$ and the sequence of indices $\bs i$.} 
For a given $\bs y$ the diagonal sequence of
polynomials determines
 the sequence of tuples $\bs y^{(i_1)}$, $\bs y^{(i_1,i_2)}$, \dots ,
$\bs y^{(i_1, \dots , i_k)}$ uniquely.

There are many diagonal sequences of polynomials associated to a given Bethe solution
and a given sequence of indices.

\section{Discrete Opers}\label{opers}
In the remaining part of the paper, 
 $\g=\g(A)$ is a simple Lie algebra of rank $r$.

Denote the coroots ${}^t\al_1^\vee, \dots ,\ {}^t\al_r^\vee
 \in \ {}^t\h$ of ${}^t\g$
by $H_1, \dots , H_r$, respectively. Let $H_1, \dots , \linebreak
H_r, \
E_1, \dots , E_r, \ F_1, \dots , F_r$ be the Chevalley generators of ${}^t\g$.
We have $[H_j,E_i] = a_{i,j} E_i$ and $[H_j,F_i] = -a_{i,j} F_i$,
where $A = (a_{i,j})$ is the Cartan matrix of $\g$.

Let $G$ be the complex simply connected
Lie group with Lie algebra $^t\g$. Let \ $B_\pm,$\ $ N_\pm$\ be 
the subgroups of $G$ with Lie algebras \ $^t\B_\pm,\ {}^t\n_\pm$,\  respectively.

\subsection{Relations in $G$}
For a non-zero complex number $u$ and $i \in \{ 1, \dots , r \}$,
consider the elements
$u^{H_i}$, $ \exp\, ( u E_i ), \ \exp\, ( u F_i)$ in $G$.
We will use the following relations. 

\begin{lemma}\label{relations}
Let $u, v$ be non-zero complex numbers. Then
\bea
u^{H_j}\ \exp \,( v \,E_i )\ = \ \exp \,( u^{a_{i,j}}v \,E_i)\ u^{H_j}\ ,
\eea
\bea
u^{H_j}\ \exp\, ( v \,F_i )\ =\ \exp \,(u^{-a_{i,j}} v \,F_i )\ u^{H_j}\ ,
\eea
\bea
\exp\, ( u \,F_i ) \ \exp\, ( v \,E_j )\  =\ \exp \,( v \,E_j )\ \exp\, ( u\, F_i )
\qquad
\text{if}  \ {} \ i \neq j \ , 
\eea
\bea
\exp\, ( u\, F_i )\ \exp \,( v\, E_i )\  =\  \exp\, ( \frac v{ 1 + uv}\, E_i )\ 
(1 + uv)^{- H_i}\  \exp\, ( \frac u{ 1 + uv} \,F_i )
\eea
if $1 + uv \neq 0$.
\hfill 
$\square$
\end{lemma}

\subsection{ D-opers }
Define the shift operator $\p$ acting on functions of $x$ by the formula
\bea
\p \ :\ g(x)\ \mapsto \ g(x + h) \ .
\eea
{\it A  discrete oper (a  d-oper)} is a difference operator of the form
$D\ =\ \partial \ - \ V$ , \
where\ $V : \C \to G$ is a rational function.

For a rational function $s : \C \to N_+$, define the action of $s$ on the
d-oper  by the formula
\bea
s \cdot D \ =\  s(x + h) \ D\ s(x)^{-1}\ =\ \p\ -\  s(x + h) \ V(x)\ s(x)^{-1}\ .
\eea
The operator $s \cdot D$ is a d-oper. The d-opers $D$ and $s \cdot D$
are called {\it gauge equivalent}.

\subsection{ Miura d-opers associated to tuples of polynomials }
In the remaining part of the paper we assume that
$\bs b = (b_{i,m})_{i,m = 1,\ i\neq m}^r$ are special symmetric parameters
given by \Ref{sym}.

Fix  dominant integral  weights  
$\bs\La = (\La_1, \dots , \La_n)$ of $\g$, complex
numbers $\bs z= (z_1, \dots , \linebreak
z_n)$. \
Introduce polynomials $T_1(x), \dots , T_r(x)$ by formula
\Ref{T}.

For given polynomials $T_1, \dots , T_r$, a
tuple $\bs y = (y_1 , \dots , y_r)$ of non-zero polynomials, and
$i \in \{1, \dots , r\}$, we 
define the rational function $R_{\bs y,  i}$  by the formula
\bea
R_{\bs y,  i} (x) \ =\ \frac{T_i(x)} { y_i(x + h )\ y_i(x)}  
\prod_{m, \ m\neq i} \  ( y_m(\, x +  b_{i,m} \,) )^{-a_{i,m}}\ .
\eea

We say that a d-oper
 $D = \p - V$ is {\it the Miura  d-oper
associated to weights $\bs\La$, numbers $\bs z$}, and the tuple
$\bs y = (y_1, \dots , y_r)$ if 
\bea
V(x) & =& \prod_{j=1}^r y_j(x + h)^{-H_j}
%\phantom{aaaaaaaaaaaaaaaaaaaaaaaaaaaaaaaaaaaaaaaaaaaaaaaaaaaaaaaaaaaaaaaaa}
\\
& \times &
 \exp \,( R_{\bs y, 1}(x)\,F_1)
\exp \,( R_{\bs y, 2}(x)\,F_2)
\cdots 
\exp \,( R_{\bs y, r}(x)\,F_r)
\prod_{j=1}^r y_j(x)^{H_j}\ .
\eea
We denote by
$D_{\bs y} = \p - V_{\bs y}$ the Miura d-oper associated to the tuple $\bs y$.

It is easy to see that if a Miura d-oper $D$ is associated to weights $\bs\La$, 
numbers $\bs z$, and a tuple $\bs y = (y_1, \dots , y_r)$ of non-zero polynomials,
then the tuple $\bs y$ is determined uniquely by $D$.

It follows easily
from Lemma \ref{relations} that the d-oper $D_{\bs y}$
does not change if the polynomials of the tuple $\bs y$ are multiplied by non-zero numbers.

For $i \in \{1, \dots , r\}$, we say that the Miura d-oper 
$D_{\bs y}$ {\it is deformable in the $i$-th direction} if there exists a non-zero
rational function $g_i : \C \to \C$ and a non-zero polynomial $\tilde y_i$ such that
\bea
\exp (g_i(x+h)\,E_i)\ D_{\bs y}\
\exp (- g_i(x)\,E_i)\ = \ D_{\bs y^{(i)}}\ ,
\eea
where $\bs y^{(i)} = ( y_1, \dots , \tilde y_i , \dots , y_r)$.
We say that $D_{\bs y}$ {\it is deformed to $D_{\bs y^{(i)}}$ with the help of
$g_i$.}

\begin{theorem}\label{main}
Let $\bs b$ be the special symmetric parameters.
Let $\bs y$ be a tuple of non-zero polynomials.
\begin{enumerate}
\item[(i)] Assume that 
equation \Ref{wron}  has a polynomial solution
$\tilde y_i$. Set  $\bs y^{(i)} = ( y_1, \dots , \tilde y_i , 
\linebreak
\dots , y_r)$. 
Then the Miura d-oper $D_{\bs y}$ is deformable in the $i$-th direction to the
Miura d-oper $D_{\bs y^{(i)}}$ with the help of
\bean\label{relation}
g_i(x) \ = \ \frac{ 1 } { y_i(x) \,\tilde y_i(x) }\
\!\!\prod_{m, \ m\neq i}   y_m( x )^{-a_{i,m}} \ .
\eean
\item[(ii)] If the tuple $\bs y$ is generic in the sense of Section
\ref{PLCP} and the Miura d-oper $D_{\bs y}$ is deformable in the $i$-th direction
to the Miura d-oper $D_{\bs y^{(i)}}$ with the help of $g_i$, then $\tilde y_i$ is
a polynomial solution of equation \Ref{wron}, and $g_i, \tilde y_i$  satisfy
\Ref{relation}.
\end{enumerate}
\end{theorem}

\begin{proof} For a scalar rational function $g_i$ we have
\bea
&&
\exp (g_i(x+h)\,E_i) \ D_{\bs y}\
\exp (- g_i(x)\,E_i)\ = \ \p - \
\\ 
&&
\phantom{aa}\ 
\prod_{j=1}^r y_j(x + h)^{-H_j} 
 \exp \,( R_{\bs y, 1}(x)\,F_1)\
\exp \,( R_{\bs y, 2}(x)\,F_2)\
\cdots 
\\
&&
\phantom{aaaaaaaaaa}\ 
\exp (\tilde  g_i(x+h)\,E_i)\ 
\exp \,( R_{\bs y, i}(x)\,F_i)\
\exp (- \tilde g_i(x)\,E_i)\ \cdots
\\
&&
\phantom{aaaaaaaaaaaaaaaaaaaaaaaaaaaaaaaa}\ 
\exp \,( R_{\bs y, r}(x)\,F_r)\
\prod_{j=1}^r y_j(x)^{H_j}\ ,
\eea
where
\bea
\tilde g_i( x )\ = \ g_i( x ) \ 
\prod_{m = 1}^r   y_m( x )^{a_{i,m}}\ .
\eea
By Lemma \ref{relations} the Miura d-oper $D$ is deformable in the $i$-th direction
only if
\bea\label{Ric}
\tilde g_i(x + h) = \tilde g_i(x) / (1 - \tilde g_i(x) R_{\bs y, i}(x)).
\eea
This equation is called {\it the $i$-th discrete Ricatti equation}, see \cite{MV4}
where the classical Ricatti equation appears in an analogous situation.

The Ricatti equation can be written as
\bean\label{ric}
&&
y_i(x + h)\ \frac {y_i(x)}
{\tilde g_i(x)}\ - \
y_i(x) \ \frac {y_i(x + h)}
{\tilde g_i(x + h)}\ =
\\
&&
\phantom{aaaaaaaaaaaaaaaaaaaaaaaaaaa}
T_i(x) \!\!
\prod_{m, \ m\neq i}   y_m(x + b_{i,m})^{-a_{i,m}}\ .
\notag
\eean
If equation \Ref{wron} has a polynomial solution $\tilde y_i$, then 
\Ref{ric} has a rational solution 
\bea
\tilde g_i (x)\ =\ \frac{ y_i(x)}{\tilde y_i(x)} .
\eea
Then $g_i$ is given by \Ref{relation}, and 
\bean\label{1}
1\ -\ \tilde g_i( x )\ R_{\bs y, i}(x)\ =\
\frac{ \tilde y_i(x + h)}{ y_i(x + h)} \
{} \frac{ y_i(x)}{\tilde y_i(x)} ,
\eean
\bean\label{2}
&&
\exp (\tilde  g_i(x+h)\,E_i)\ 
\exp \,( R_{\bs y, i}(x)\,F_i)\
\exp (- \tilde g_i(x)\,E_i)\ = \
\\
&&
\left( \frac{ y_i(x + h)} { \tilde y_i(x + h)}
\right)^{H_i} \left( \frac{\tilde y_i(x)} { y_i(x)}\right)^{H_i}\times
\notag
\\
&&
\exp \,\left( \frac {\tilde y_i (x)\ T_i(x)}
{\tilde y_i (x + h)\ (y_i (x))^2}
\prod_{m, \ m\neq i}   y_m( x + b_{i,m})^{-a_{i,m}}\,F_i\right) =
\notag
\\
&&
\left( \frac{ y_i(x + h)} { \tilde y_i(x + h)}
\right)^{H_i}
\exp \,( R_{\bs y^{(i)}, i}(x)\,F_i)\
\left( \frac{\tilde y_i(x)} { y_i(x)}\right)^{H_i} .
\notag
\eean
Using the last formula and Lemma \ref{relations}
 we easily conclude that the Miura d-oper $D_{\bs y}$ is deformed in the $i$-th
direction to the Miura d-oper $D_{\bs y^{(i)}}$ with the help of $g_i$ given by
\Ref{relation} if $\tilde y_i$ is a polynomial solution of \Ref{wron}. 
This proves part (i) of the theorem.

To prove part (ii) write \Ref{ric} as
\bea
\frac 1
{\tilde g_i(x)}\ - \
 \frac 1
{\tilde g_i(x + h)}\ =
\frac{T_i(x)}
{y_i(x + h)\ y_i(x)}
\prod_{m, \ m\neq i}   y_m( x + b_{i,m})^{-a_{i,m}}\ .
\eea
Let $\tilde g_i(x)$ be a rational solution of this equation. Since $\bs y$ is generic,
the poles of $1/\tilde g_i(x)$ are located at zeros of $y_i(x)$ and all poles are simple.
Hence $\tilde y_i(x) = y_i(x)/ \tilde g_i(x)$ 
is a polynomials solution of \Ref{wron}. 
Then formulas \Ref{1}, \Ref{2} hold and part (ii) is proved.
\end{proof}

\begin{corollary}
Let the  Miura d-oper
$D_{\bs y}$ be associated to
weights $\bs\La$, numbers $\bs z$, and the tuple $\bs y = (y_1, \dots , y_r)$.
Assume that the  tuple $\bs y = (y_1, \dots , y_r)$ is generic in the sense of
Section \ref{PLCP}. Then
 $D_{\bs y}$ is deformable in all directions from 1 to $r$
if and only if the tuple $\bs y$
represents  a Bethe solution of the Bethe ansatz equations 
associated to \ $\bs z,\ \bs \La,\ \La_\infty = \sum_{i=1}^n \La_i - 
\sum_{i=1}^r l_i\al_i$,\ $l_i = \deg y_i$, and special symmetric parameters $\bs b$.
\end{corollary}

Let the  Miura d-oper
$D_{\bs y}$ be associated to
weights $\bs\La$, numbers $\bs z$, and the tuple $\bs y = (y_1, \dots , y_r)$.
Let the  tuple $\bs y = (y_1, \dots , y_r)$ 
represent   a Bethe solution of the Bethe ansatz equations 
associated to  $\bs z, \bs \La, \La_\infty$, and special symmetric parameters
$\bs b$.  
Let $\om^0 $ be the variety 
of all Miura d-opers each of 
which can be obtained from $D_{\bs y}$ by a sequence of deformations
in directions $i_1, \dots , i_k$ where $k$ is a positive integer and all
$i_j$ lie in $\{1, \dots , r\}$.

\begin{corollary}\label{om0}
The variety $\om^0$ is isomorphic to the population $P_{\bs y}$ 
of solutions of Bethe ansatz equations, where $P_{\bs y}$ is the population
originated at $\bs y$.
\end{corollary}

\section{Solutions of Difference Equations}\label{solutions}

Let $D_{\bs y} = \p - V_{\bs y}$ be the Miura d-oper associated to a Bethe solution
$\bs y$ of the Bethe ansatz equations associated to
special symmetric parameters $\bs b$. Let $P_{\bs y}$ be the 
population of solutions originated at $\bs y$.

 In this section
we prove that the difference equation 
\bean\label{eqn1}
Y(x+h)\ =\ V_{\bs y}(x)\ Y(x)
\eean
has a $G$-valued rational solution. We will write that solution 
explicitly in terms of coordinates of tuples composing the population.

Note that if $Y(x)$ is a solution and $g \in G$, then $Y(x) g$ is a solution too.

First we give a formula for a solution of equation \Ref{eqn1}
for d-opers associated to $\g$ of type $A_r$,
and then we consider more general formulas for solutions
which do not use the structure of the Lie algebra.

Let $Y$ be a solution of equation \Ref{eqn1}. Define 
\bea
\bar Y(x)\ =\ \prod_{j=1}^r y_j(x)^{H_j}\ Y(x) \ .
\eea
Then $\bar Y$ is a solution of the equation
\bean\label{eqn}
\bar Y(x+h)\ =\ \bar V_{\bs y}(x)\ \bar Y(x)
\eean
where
\bea 
\bar V_{\bs y}(x) =
 \exp \,( R_{\bs y, 1}(x)\,F_1)\,
\exp \,( R_{\bs y, 2}(x)\,F_2)
\cdots 
\exp \,( R_{\bs y, r}(x)\,F_r) .
\eea

\subsection{The $A_r$ d-opers and solutions of Bethe ansatz equations}\label{AO} 
In this section let $\g = sl_{r+1}$ be the Lie algebra of type $A_r$.
Then ${}^t\g = sl_{r+1}$. We have $(\al_i, \al_i) = 2$ for all $i$.
We fix the order of simple roots of $sl_{r+1}$ such that
\bea
(\al_1, \al_2) = (\al_2, \al_3) = \dots = (\al_{r-1}, \al_{r}) = - 1 .
\eea
We start with two examples.

Let $\g = sl_2$. Let $\bs y = (y_1)$ represent a Bethe solution of the $sl_2$ 
Bethe ansatz equations associated to  $\bs z,\ \bs \La,\ \La_\infty$.
Let $y_1^{(1)}$ be the diagonal sequence
of polynomials associated to $\bs y$ and the sequence of indices $(1)$, \ 
in other words,
\bea
y_1(x+h)\ y_1^{(1)}(x)\ -\ y_1(x)\ y_1^{(1)}(x+h)\ 
 =\ T_1(x)\ .
\eea
 Then
\bea
\bar Y \ =\ \exp\,( \frac{y_1^{(1)}}{ y_1} F_1) 
\eea
is a solution of the difference equation  \Ref{eqn} with values in $\mathrm{SL}\,(2,\C)$.
Indeed,
\bea
&&
(\p - \exp(\, \frac{ T_1(x)}{y_1(x) y_1(x+h)}\,F_1 ) )\ \exp\,( \frac{y_1^{(1)}(x)}
{ y_1(x)} F_1) \ =
\ \exp\,( \frac{y_1^{(1)}(x+h)}{ y_1(x+h)} F_1) \times
\\
&&
  (\p - \exp(\, (\,
 - \frac{y_1^{(1)}(x+h)}{ y_1(x+h)} + \frac{y_1^{(1)}(x)}{ y_1(x)}
+  \frac{T_1(x)}{y_1(x) y_1(x+h)})\,F_1\,)\,)\ =
\\
&&
\exp\,( \frac{y_1^{(1)}(x+h)}{ y_1(x+h)} F_1) \
( \p - \id )\ .
\eea

Let $\g = sl_3$. Let $\bs y = (y_1, y_2)$ represent a Bethe solution of the Bethe ansatz
equations   associated to  $\bs z, \bs \La, \La_\infty$, and
special symmetric parameters $\bs b$.
Let $ y_1^{(1)},  y_2^{(1,2)}$ be the diagonal sequence
of polynomials associated to $\bs y$  and the sequence of indices $(1,2)$, \ 
in other words,
\bea
 y_1(x+h)\ y_1^{(1)}(x)\ -\  y_1(x)\ y_1^{(1)}(x+h)\ &=&\ T_1(x)\ y_2(x+h)\ ,
\\
 y_2(x+h) \ y_2^{(1,2)}(x)\ -\  y_2(x)\ y_2^{(1,2)}(x+h)\ &=&\ T_2(x)\ 
y_1^{(1)}(x)\ .
\eea
Let $y_2^{(2)}$ be the  diagonal sequence
of polynomials associated to $\bs y$  and the sequence of indices $(2)$, \
in other words,
\bea
 y_2(x+h)\ y_2^{(2)}(x) \ - \ y_2(x) \ y_2^{(2)}(x+h) \ =\ T_2(x)\ y_1(x)\  .
\eea
Then
\bea
\bar Y(x)\ =\ \exp\,( \frac{y_1^{(1)}(x)}{ y_1(x)} F_1)\
\exp\,( \frac{ y_2^{(1,2)}(x)}{ y_2(x)} [F_2,F_1])\
\exp\,( \frac{ y_2^{(2)}(x)}{ y_2(x)} F_2 )
\eea
is a solution of the difference equation \Ref{eqn} with values in $\mathrm{SL}\,(3,\C)$.
Indeed, we have
\bea
&&
(\p - \bar V_{\bs y}(x))\ \exp\,( \frac{y_1^{(1)}(x)}{ y_1(x)} F_1)\ =\
\\
&&
\exp( \frac{y_1^{(1)}(x+h)}{ y_1(x+h)} F_1)(\p  -
\exp( (
 - \frac{y_1^{(1)}(x+h)}{ y_1(x+h)} + \frac{y_1^{(1)}(x)}{ y_1(x)}
+  \frac{T_1(x)y_2(x+h)}{y_1(x) y_1(x+h)}) F_1 ) \times
\\
&&
\phantom{aaaaaaaaaaaa}
\exp\,( \frac{ T_2(x) y_1^{(1)}(x)}{ y_2(x) y_2(x+h)} [F_2,F_1])\
\exp\,( \frac{ T_2(x) y_1(x)}{ y_2(x) y_2(x+h)} F_2)\ ) =
\\
&&
\exp( \frac{y_1^{(1)}(x+h)}{ y_1(x+h)} F_1) 
\times
\\
&&
\phantom{aaaaaaaa}
(\p  -
\exp( \frac{ T_2(x) y_1^{(1)}(x)}{ y_2(x) y_2(x+h)} [F_2,F_1])
\exp( \frac{ T_2(x) y_1(x)}{ y_2(x) y_2(x+h)} F_2)) ,
\eea
\bea 
&&
(\p  -
\exp( \frac{ T_2(x) y_1^{(1)}(x)}{ y_2(x) y_2(x+h)} [F_2,F_1])
\exp( \frac{ T_2(x) y_1(x)}{ y_2(x) y_2(x+h)} F_2)) \times
\\
&&
\phantom{aaaaaaaaaaaaaaaaaaaaaaaaaaaaaaaaaaaaaaaa}
\exp\,( \frac{ y_2^{(1,2)}(x)}{ y_2(x)} [F_2,F_1]) =
\\
&&
\phantom{aaaaaaaaaaaaa}
\exp\,( \frac{ y_2^{(1,2)}(x+h)}{ y_2(x+h)} [F_2,F_1]) 
\ (\p  -
\exp( \frac{ T_2(x) y_1(x)}{ y_2(x) y_2(x+h)} F_2)) ,
\eea
and
\bea
&&
(\p  - \exp( \frac{ T_2(x) y_1(x)}{ y_2(x) y_2(x+h)} F_2))\
\exp\,( \frac{ y_2^{(2)}(x)}{ y_2(x)} F_2 ) =
\\
&&
\phantom{aaaaaaaaaaaaaaaaaaaaaaaaaaaaaaaa}
\exp\,( \frac{ y_2^{(2)}(x+h)}{ y_2(x+h)} F_2 )
\ (\p  - \id )\  .
\eea

Consider the general case.
Let $\g = sl_{r+1}$.
Let $\bs y = (y_1, \dots , y_r)$ represent a Bethe solution of the Bethe ansatz
equations
associated to  $\bs z, \bs \La, \La_\infty$, and special symmetric parameters $\bs b$.
For $i = 1, \dots , r$,\ let 
$ y_i^{(i)},\  y_{i+1}^{(i,i+1)},\ \dots ,\
 y_{r}^{(i,\dots , r)}$  be the diagonal sequence
of polynomials associated to $\bs y$  and the sequence of indices $(i,i+1, \dots , r)$, \ 
in other words,
\bea
 y_i(x+h)\ y_i^{(i)}(x)\ -\  y_i(x)\ y_i^{(i)}(x+h)
 \ =\ T_i(x)\ y_{i-1}(x) \  y_{i+1}(x+h)\  ,  
\eea
\bea
 y_{i+1}(x+h) y_{i+1}^{(i, i+1)}(x) -  y_{i+1}(x) y_{i+1}^{(i, i+1)}(x+h)
  = T_{i+1}(x)\ y^{(i)}_{i}(x)   y_{i+2}(x+h)  ,\ {} \ \dots ,  
\eea
\bea
y_{r-1}(x+h)y_{r-1}^{(i,\dots , r-1)}(x)-
y_{r-1}(x)y_{r-1}^{(i,\dots , r-1)}(x+h) = T_{r-1}(x)\
y_{r-2}^{(i, \dots , r-2)}(x) y_{r}(x+h),
\eea
\bea
y_{r}(x+h)\ y_{r}^{(i,\dots , r)}(x) -
y_{r}(x)\ y_{r}^{(i,\dots , r )}(x+h)
\ =\ T_r(x)\ y_{r-1}^{(i, \dots , r-1)}(x) .
\eea
Define $r$ functions $ Y_1, \dots , Y_r$
of $x$ with values in $\mathrm{SL}\,(r+1,\C)$ by the formulas
\bea
Y_i(x)\ = \ \prod_{j=i}^r \ \exp\,(\, \frac { 
y_j^{(i,\dots ,j)}(x)}{ y_j(x)}
\ [F_j,[F_{j-1},[...,[F_{i+1},F_i]...]]]\,)\ .
\eea

Note that inside each product the factors commute.

\begin{theorem}\label{sl}
The product $Y_1 \cdots Y_r$ is a solution of the difference equation
\Ref{eqn} with values in $\mathrm{SL}\,(r+1,\C)$.
\end{theorem}

The proof is straightforward. One shows that
\bea
&&
( \p\ - \ 
\exp \,( R_{\bs y, i}(x)\,F_i)
\cdots 
\exp \,( R_{\bs y, r}(x)\,F_r)) \ Y_i(x)\ =
\\
&&
\phantom{aaaaaa}
Y_i(x+h)\ 
( \p\ - \ 
\exp \,( R_{\bs y, i+1}(x)\,F_{i+1})
\cdots 
\exp \,( R_{\bs y, r}(x)\,F_r)) \ .
\eea

\subsection{General formulas for solutions} 
Let $U$ be a complex finite dimensional representation of $G$. Let $u_{low}$ be 
a lowest weight vector of $U$, \ {} ${}^t\n_-\, u_{\mathrm{low}} = 0$.

Let $\bs y = (y_1, \dots , y_r)$ represent a Bethe solutions of the Bethe  ansatz
equations
associated to  $\bs z, \bs \La, \La_\infty$, and special symmetric parameters $\bs b$.
We solve the difference equation \Ref{eqn1} with values in $U$.

Let $\bs i = (i_1, i_2, \dots , i_k), \ 1 \leq i_j \leq r,$ 
be a sequence of natural numbers. Let
$\bs y^{(i_1)} = ( y_1^{(i_1)}, \dots , y_r^{(i_1)})$,
$\bs y^{(i_1, i_2)} = ( y_1^{(i_1,i_2)}, \dots , y_r^{(i_1,i_2)})$,
\dots ,
$\bs y^{(i_1,  \dots , i_k)} = ( y_1^{(i_1, \dots , i_k)}, \dots , 
y_r^{(i_1, \dots , i_k)})$ be a sequence of tuples
associated to the Bethe solution  $\bs y$
and the sequence of indices $\bs i$.

\begin{theorem}\label{Solutions}
The $U$-valued function
\bea\label{Sol}
Y(x)& =&
\exp\,\left( -
\frac{ 1 } { y_{i_1}(x) \, y^{(i_1)}_{i_1}(x) }
\prod_{m, \, m\neq i_1}  ( y_m(x) )^{-a_{i_1,m}} E_{i_1} \right)
\\
& \times &
\exp\,\left( -
\frac{ 1 } { y^{(i_1)}_{i_2}(x) \, 
y^{(i_1,i_2)}_{i_2}(x) }\
\prod_{m, \, m\neq i_2}   ( y^{(i_1)}_m(x) )^{-a_{i_2,m}}  E_{i_2} \right)\ \cdots
\\
&\times &
\exp\,\left( -
\frac{ 1 } { y^{(i_1,\dots,i_{k-1})}_{i_k}(x) 
y^{(i_1,\dots,i_k)}_{i_k}(x) }
\prod_{m, \, m\neq i_k} ( y^{(i_1,\dots, i_{k-1})}_m(x) )^{-a_{i_k,m}}
 E_{i_k} \right)
\\
&\times &
\prod_{j=1}^r\ ( y^{(i_1,\dots,i_k)}_j(x) )^{- H_j}\ u_{\mathrm{low}} \ .
\eea
is a solution of the difference equation \Ref{eqn1}.
\end{theorem}

The proof is straightforward and follows from Theorem \ref{main}.

\begin{corollary}\label{Rationality}
Every coordinate of every solution of the difference equation \Ref{eqn1} with values in a
finite dimensional representation of \ $G$ can be written as a rational function
$R(f_1, \dots , f_N)$ of 
suitable polynomials $f_1, \dots , f_N$ which appear as coordinates
of tuples in the $\g$ population $P_{\bs y}$ generated at $\bs y$.
\end{corollary}

Since $G$ has a faithful finite dimensional representation, the solutions of the difference
equation \Ref{eqn1} with values in $G$ also can be written as rational functions
of coordinates of tuples of $P_{\bs y}$, cf. Section \ref{AO}.

\begin{corollary}\label{rationality}
Let $\bs y = (y_1, \dots , y_r)$ represent a Bethe solution of the Bethe  ansatz
equations
associated to  $\bs z, \bs \La, \La_\infty$, and special symmetric parameters $\bs b$.
Then there
exists a $G$-valued rational function $Y : \C \to G$ satisfying equation
\Ref{eqn1}.
\end{corollary}

\section{Miura d-opers and flag varieties}

\subsection{Theorem on isomorphism}

Let $\bs y = (y_1, \dots , y_r)$ represent a Bethe solution of the Bethe ansatz equations
associated to  $\bs z, \bs \La, \La_\infty$,
 and special symmetric parameters $\bs b$.
Let $D_{\bs y} = \p - V_{\bs y}$ be the Miura d-oper associated to $\bs y$.
Consider the variety $\om$ of all Miura d-opers gauge equivalent to $D_{\bs y}$. 
If $D' \in \om$, then there exists a rational function
$v : \C \to N_+$ such that $D' \ =\ v(x+h)\, D_{\bs y}\, v(x)^{-1}$. 
In that case we denote $D'$ by $D^v$. 

The variety of pairs 
\bea
\widehat {\om} \ =\ 
\{ (D^v, v)\ | \ D^v \in \om \}
\eea
will be called {\it the variety of marked Miura d-opers gauge equivalent to } $D_{\bs y}$.
We have the natural projection $\pi : \widehat {\om} \to \om,\ (D^v,v) \mapsto D^v$.
We will show below that $\pi$ is an isomorphism.

Let $\om^0 \subseteq \om$ be the subvariety 
of all Miura d-opers each of which can be 
obtained from $D_{\bs y}$ by a sequence of deformations
in directions $i_1, \dots , i_k$ where $k$ is a non-negative integer and all
$i_j$ lie in $\{1, \dots , r\}$. By Corollary \ref{om0} the subvariety
$\om^0$ is isomorphic to the population of solutions of Bethe ansatz
equations
associated to special symmetric parameters $\bs b$ and originated at $\bs y$.
We will show below that $\om^0 = \om$.

Assume that 
$D' \in \om^0$ and $D'$ is obtained from $D_{\bs y}$ by a sequence of deformations
in directions $i_1, \dots , i_k$ where $k$ is a non-negative integer and all
$i_j$ lie in $\{1, \dots , r\}$. Then there exist
scalar rational functions $g_1, \dots , g_k$ with the following properties.
For $j = 1, \dots , k$, define a rational $N_+$-valued function $v_j : \C \to N_+$,
\bean\label{presentation}
v_j(x) \ =\  \exp\, (g_j(x) E_{i_j}) \cdots \exp\, (g_2(x) E_{i_2}) \exp\, (g_1(x) E_{i_1}) \ .
\eean
Then $D^{v_j} \in \om^0$ and $D' = D^{v_k}$.
The set of all pairs $(D^{v_k}, v_k)$ such that $k $ is a non-negative integer,
$v_k$ is given by the above construction, and
$D^{v_k} \in \om^0$, will be called {\it the variety of specially
marked Miura d-opers gauge equivalent to } $D_{\bs y}$
and denoted by $\widehat{\om^0}$.
Clearly we have $\widehat {\om^0} \subseteq \widehat{\om}$.

\bigskip

Let $\Bbb P^1$ be the complex projective line. Consider $D_{\bs y}$ as a
discrete connection $\nabla_{\bs y}$ on the trivial principal $G$-bundle 
$p : G \times \Bbb P^1 \to \Bbb P^1$.  Namely, by definition
a section
\bea
U  \ \to\ G \times U ,
\qquad
x\ \mapsto\ Y(x) \times x ,
\eea 
of $p$ over a  subset $U \subset \C \subset \bP$ 
is called {\it horizontal} if the $G$-valued function
$Y(x)$ is a solution of the difference equation \Ref{eqn1},
$Y(x+h)\, = \,V_{\bs y}(x) \,Y(x)$.  
%the horizontal sections of the discrete connection
%$\nabla_{\bs y}$ are $G$-valued meromorphic
%solutions of the difference equation \Ref{eqn1},
%$Y(x+h)\, = \,V_{\bs y}(x) \,Y(x)$.  
By Corollary \ref{rationality} equation \Ref{eqn1} has a rational solution 
$Y(x)$. For any $g \in G$ the rational $G$-valued function
$Y(x)g$ is a solution of the same equation  too. 
A point $x_0 \in \C$ will be called {\it regular} if $x_0$ is a regular 
point of the rational functions $Y(x)$ and $V_{\bs y}(x)$.

Let $x_0\in \C$ be a regular point. Let
$g$ be an element of $G$. Then $\nabla_{\bs y}$ has a rational horizontal section $s$ 
such that $s(x_0)= g$.

It is easy to see that if the values of two rational horizontal sections are equal
at one point, then the sections are equal.

\bigskip

Consider the trivial bundle  $p' : (G/B_-) \times \Bbb P^1 \to  \bP$
associated to the bundle $p$. The fiber of $p'$ is the flag variety $G/B_-$. 
The discrete connection $\nabla_{\bs y}$ induces
a discrete connection $\nabla_{\bs y}'$ on $p'$. 

The variety $\Gamma$ of rational horizontal sections of the discrete connection
$\nabla_{\bs y}'$ is identified with the fiber $(p')^{-1}(x_0)$ over any 
regular point $x_0$.  Thus, $\Gamma$ is isomorphic to $G/B_-$.

Any $G$-valued rational function $v$ defines a section
\bean\label{section}
S_v \ :\ x\ \mapsto\ v(x)^{-1} B_- \times x
\eean
of $p'$ over the set of regular points of $v$.
The section $S_v$ is also well defined over the poles of $v$ since $G/B_-$ 
is a projective variety.

If $D^v \in \om$, then the section $S_v$ is horizontal with respect to
$\nabla_{\bs y}'$. This follows from the fact that the function $V_{\bs y}$
takes values in $B_-$.
Thus we have a map 
\bea
S \ {} :\ {} \widehat \om \to \Gamma , 
\qquad
(D^v, v) \mapsto S_v \ .
\eea

\begin{theorem}\label{isomorphism}
The map $ S : \widehat \om \to \Gamma $ is an isomorphism and 
$\widehat{\om^0} = \widehat \om$.
\end{theorem}

\begin{proof}
Let $(D^{v_1}, v_1), (D^{v_2}, v_2) \in \widehat \om$. 
Assume that the images of 
$(D^{v_1}, v_1)$ and $ (D^{v_2}, v_2)$  under the map $S$
coincide. Assume that $v_1, v_2, V_{\bs y}$ are regular
at $x_0\in \C$. The equality $S_{v_1}(x_0)= S_{v_2}(x_0)$
means that  $v_1(x_0)^{-1} B_- = v_2(x_0)^{-1} B_-$. Then
$v_1(x_0) = v_2(x_0)$. Hence $v_1 = v_2$ and $D^{v_1} = D^{v_2}$.
That proves the injectivity of $S$. 

Let $x_0$ be a regular point of $V_{\bs y}$ in $\C$. For any $u \in  N_+$
there exists a rational \linebreak
$N_+$-valued function $v$ such that $v(x_0) = u$,
 $D^v \in \om^0$.
Indeed, every $u \in N_+$ is a product of elements of the form 
$e^{c_i E_i }$ for $i\in\{1, \dots , r\}$ and $c_i\in\C$.
Every $c_i$ can be taken as the initial condition for a solution of the suitable
$i$-th discrete Ricatti equation. 

Thus the set
\bea
Im (x_0) \  =\ \{ S_v(x_0) \in \ (G/B_-)\times x_0 \ | \ (D^v, v) \in \widehat{\om^0} \}
\eea
contains the set $((N_+ B_-)/ B_- )\times x_0\ \subset\ (G/B_-)\times x_0$.
It is easy to see that the set $Im(x_0)$ is closed in
$(G/B_-)\times x_0$ as the image of $\widehat{\om^0}$ with respect to $S$.
 On the other hand
the set $((N_+ B_- )/B_-)\times x_0$ is dense in $(G/B_-)\times x_0$. Hence
$Im(x_0) = (G/B_-)\times x_0$, and $\widehat{\om^0} = \widehat\om$ 
since the map $S$ is injective.
\end{proof}

\subsection{Remarks on the isomorphism}\label{REMARKS}
Let $\g$ be a simple Lie algebra. 
Let $\bs y^0$ be a Bethe solution of the Bethe ansatz equations
associated to special symmetric parameters $\bs b$.
Theorem \ref{isomorphism} says that the variety
$\widehat\tN$ is isomorphic to the flag variety $G/B_-$. 
Here are some comments on that isomorphism.

The isomorphism is constructed in two steps. 
If $(D^v, v) \in \widehat\tN$ is a marked Miura d-oper, then we assign to it
the  section $S_v \in \Gamma$ by formula \Ref{section}.
We choose a regular point $x_0\in \C$,
and assign to the section $S_v$ its value $S_v(x_0) \in (G/B_-) \times x_0$ at $x_0$.
The resulting composition 
\bea
\phi_{\bs y^0, x_0}\ :\ \widehat\tN \ \to \ G/B_- 
\eea
is an isomorphism according to Theorem \ref{isomorphism}.

\begin{lemma}\label{dep 1}
If $x_0, x_1 \in \C$ are regular points,
then there exists an element $g\in B_-$ such that \ {}\
$\phi_{\bs y^0, x_1}\ = \ g \ \phi_{\bs y^0,  x_0}$.
\end{lemma}

\begin{proof} Let $Y$ be the $G$-valued rational
solution of equation
\Ref{eqn1} such that $Y(x_0) = \id$. Then $Y(x) \in B_-$ 
for all  $x$. If $(D^v, v)  \in \widehat\tN$,
then $S_v$ is a horizontal section
of $\nabla_{\bs y^0}'$. Thus it has the form
$ x \ \mapsto\ (Y(x) u  B_-)\times x$ for a suitable element $u\in G$.
Hence $\phi_{\bs y^0, x_0} ( \bs y ) \ =\  Y(x_0) u  B_-$
and
$\phi_{\bs y^0, x_1} ( \bs y )\ =\ Y(x_1) u B_-$.
We conclude that
$\phi_{\bs y^0, x_1}\ =\ Y(x_1)Y(x_0)^{-1} \,\phi_{\bs y^0, x_0}$.
\end{proof}

\section{Bruhat Cells}\label{cells}

\subsection{Properties of Bruhat cells} Let $\g$ be a simple Lie algebra. 
For an element $w$ of the Weyl group $W$, the set
\bea
B_w\ = \ B_-  w B_- \ {} \subset \ {} G/B_-
\eea
is called {\it the Bruhat cell} associated to $w$. The Bruhat cells
form a cell decomposition of the flag variety $G/B_-$.

For $w\in W$ denote $l(w)$ the length of $w$. We have 
$\dim B_w = l(w)$.

Let $s_1, \dots , s_r \in W$ be the generating reflections of the Weyl group.
For $v \in G/B_-$ and $i \in \{1, \dots , r\}$ consider the rational
curve 
\bea 
\C\ \to\ G/B_-, 
\qquad
 c \ \mapsto \  e^{ c E_i}v \ .
\eea
 The limit of
$e^{ c E_i}\,v$ is well defined as $c \to \infty$, since
$G/B_-$ is a projective variety.

We  need the following standard property of Bruhat cells.

\begin{lemma}\label{BRUHAT}
Let $s_i, w \in W$ be such that $l( s_i w ) = l(w) + 1$. Then
\bea
B_{ s_i w }\ =\
\{\ e^{ c E_i }v\ | \ v \in B_w,\ c\in \{\bP - 0\}\ \} \ .
\eea
\hfill
$\square$
\end{lemma}

\begin{corollary}\label{Bruhat}
Let $w = s_{i_1} \cdots s_{i_k}$ be a reduced decomposition of $w\in W$. Then
\bea
B_{ w } =
\{ \lim_{c_1 \to c^0_1}\,\dots\,\lim_{c_k \to c^0_k}
e^{ c_1 E_{i_1} }\,\cdots \,e^{ c_k E_{i_k} }B_- \ \in \ G/B_- \  | \  c^0_1, \dots , 
c^0_k
\in \{\bP - 0\} \}.
\eea
\end{corollary}

\begin{corollary}\label{Bruhat-2}
Let $s_{i_1} \cdots s_{i_k}$ be an element in $ W$.
Let $c^0_1, \dots , c^0_k \in \bP$.  Then the element
\bea
 \lim_{c_1 \to c^0_1}\,\dots\,\lim_{c_k \to c^0_k}
e^{ c_1 E_{i_1} }\,\cdots \,e^{ c_k E_{i_k} }B_- \ \in \ G/B_- 
\eea
belongs to the union of the Bruhat cells $B_w$ with $l(w)\leq k$.
\end{corollary}

\subsection{Populations and Bruhat cells}
Let $P$ be a population of solutions of the
Bethe ansatz equations associated to  integral dominant
weights $\bs\La$, numbers $\bs z$, special symmetric parameters $\bs b$.

%Let $T_1, \dots , T_r$ be the polynomials defined by \Ref{T}.
Let $\bs y^0 = (y^0_1, \dots , y^0_r) \in P$ 
be a point of the population with
 $l_i = \deg\, y^0_i$ for  $i = 1, \dots , r$.
Assume that the weight at infinity of $\bs y^0$,
\bea
\La_\infty \ =\ \sum_{i=1}^n \La_i\ - \
\sum_{i=1}^r\ l_i\,\al_i\ ,
\eea
is integral dominant, see Section \ref{crit pts}. Such $\bs y^0$ exists according to
\cite{MV1}, \cite{MV2}.
For $w\in W$ consider the weight $w \cdot \La_\infty$, where
$w \cdot $ is the shifted action of $w$ on $\h^*$. Write
\bea
w \cdot \La_\infty \ =\ \sum_{i=1}^n \La_i\ - \
\sum_{i=1}^r\ l^w_i\,\al_i .
\eea
Set
\bea
P_w\ =\ \{\ \bs y = (y_1, \dots , y_r)\,\in \,P\ | \ 
\deg y_i = l^w_i, \ i = 1, \dots , r \ \}\ .
\eea
Clearly, $P = \cup_{w\in W} P_w$, and $P_{w_1} \cap P_{w_2} = \emptyset$
if $w_1\neq w_2$.

Consider the trivial bundle  $p' : (G/B_-) \times \bP \to  \bP$ 
with  the discrete
connection $\nabla_{\bs y^0}'$. Consider the Bruhat cell decomposition
of fibers of $p'$. 

Assume that $x_0\in \C$ is a regular point of the Miura d-oper
$D_{\bs y^0} = \p - V_{\bs y^0}$ and $x_0$ is
a regular point of the $G$-valued rational solutions of the associated
difference equation \Ref{eqn1}, $Y(x+h)\, =\, V_{\bs y^0}(x)\, Y(x)$.
Let 
\bea
\phi_{\bs y^0, x_0}\ :\ \widehat\tN\ \to \ G/B_-
\eea
be the isomorphism defined in Section \ref{REMARKS}.
Let 
\bea
\pi : \widehat\tN \to \tN , 
\qquad
(D^v , v) \mapsto D^v ,
\eea
be the natural projection.
Let 
\bea
\xi : \tN \to P , 
\qquad
D_{\bs y} \mapsto \bs y ,
\eea
be the isomorphism of Corollary \ref{om0}.

\begin{theorem}\label{second main}
For every $w \in W$, the composition 
$ \xi \pi  \phi_{\bs y^0, x_0}^{-1} : G/B_- \to P$, restricted to 
the Bruhat cell $B_{w^{- 1}}$, is a 1-1 epimorphism of $B_{w^{- 1}}$ onto
$P_w$.
\end{theorem}

\begin{corollary}\label{proj iso}
The projection $\pi : \widehat\tN \to \tN$ is an isomorphism, i.e.
if $(D^{v_1}, v_1),\linebreak
 (D^{v_2}, v_2) \in \widehat\tN$ are such that
$D^{v_1} = D^{v_2}$, then $v_1 = v_2$.
\end{corollary}

\begin{corollary}\label{ISO}
Let $P$ be a population of
solutions of the Bethe ansatz equations associated to
integral dominant $\g$-weights $\La_1, \dots , \La_n$, complex numbers
$z_1, \dots , z_n$, special symmetric parameters $\bs b$. Then $P$ is isomorphic to
the flag variety $G/B_-$ of the Langlands dual algebra ${}^t\g$.
\end{corollary}

\begin{corollary}\label{main corollary}
Let $\La_1, \dots , \La_n, \La_\infty$ be integral dominant $\g$-weights.
Let $z_1, \dots , z_n$ be complex numbers. Let $w\in W$. 
Consider the Bethe ansatz equations associated to
$\La_1, \dots , \La_n, w\cdot\La_\infty = \sum_{i=1}^n\La_i - \sum_{i=1}^rl_i^w\al_i$, 
$z_1, \dots , z_n$,
and special symmetric parameters $\bs b$. 
A solution of the Bethe ansatz equations is a collection of
complex numbers $\bs t = (t^{(i)}_j),\ i = 1, \dots , r,\
j = 1, \dots , l^w_i$. 
Let $K$ be a connected component of the set of solutions of the Bethe ansatz
equations.
For each $\bs t \in K$ consider the tuple $\bs y_{\bs t} \in (\C[x])^r$ 
of monic polynomials
representing the solution $\bs t$. Then the closure of the
set $\{\ \bs y_{\bs t} \ | \ \bs t \in K \ \}$ in $(\C[x])^r$ 
is an $l(w)$-dimensional cell.
\end{corollary}

\subsection{Proof of Theorem \ref{second main}}

\begin{lemma}\label{invariance}
For $w \in W$,  the subset $B_w \times \bP \ \subset \
(G/B_-) \times \bP$ is invariant with respect to
the discrete connection $\nabla_{\bs y^0}'$.
\end{lemma}

\begin{proof} Let $Y$ be the rational $G$-valued solution of the equation
$D_{\bs y^0} Y = 0$ such that $Y(x_0) = \id$.  Then $Y(x) \in B_-$ 
for all  $x$. The rational horizontal sections of 
$\nabla_{\bs y^0}'$ have the form
$ x \ \mapsto\ (Y(x)\,u  B_-)\times x$ for a suitable element $u\in G$.
If $u B_- \ \in B_w$, then
$Y(x)u B_-\  \in B_w$ for all $x$.
\end{proof}

Let $w = s_{i_k} \cdots s_{i_1}$ be a reduced decomposition of $w\in W$. 
For $d = 1, \dots , k$ set
\bea
(s_{i_d} \cdots s_{i_1}) \cdot \La_\infty \ =\
\sum_{i=1}^n \La_i\ - \ \sum_{i=1}^r\ l^d_i\,\al_i .
\eea
From \cite{BGG} it follows that $l^{1}_{i_1} > l_{i_1}$ and
$l^{d}_{i_d} > l^{d-1}_{i_d}$ for $d = 2, \dots , k$.

Let $\bs i = (i_1, \dots , i_k), \ 1 \leq i_j \leq r$, be a sequence of integers.
We consider the map 
$Y_{\bs y^0, \bs i} \ :\ (\Bbb P^1)^k \ \to \PCr$ introduced in Section \ref
{general procedure} for special symmetric parameters $\bs b = (b_{i,j})$.
Its image is denoted by $P_{\bs y^0, \bs i}$.
The image of a point $(c_1, \dots , c_k) \in (\Bbb P^1)^k$ under this map is denoted
by $\bs y^{k;\, c_1, \dots , c_k}$. We repeat the definition of 
$\bs y^{k;\, c_1, \dots , c_k}$ in terms convenient for our present purposes.

We assume that the tuple $\bs y^0$ is a tuple of monic polynomials.
For $d = 1, \dots , k$ we define by induction on $d$
a family of tuples of polynomials  depending on
parameters $c_1, \dots , c_d \in \Bbb P^1$. 
Namely, let $\tilde y_{i_1}$ be a polynomial satisfying equation
\bea
y^0_{i_1}(x+h)\, \tilde y_{i_1}(x)\ -\
y^0_{i_1}(x) \,\tilde y_{i_1}(x+h)
\ = \ T_{i_1}(x)\! \prod_{j,\ j\neq i_1} ( \,y^0_j(x+b_{i_1,j})\, )^{- a_{i_1,j}}\ .
\eea
We fix $\tilde y_{i_1}$ assuming that the coefficient of $x^{l_{i_1}}$ in $\tilde y_{i_1}$
is equal to zero.
Set
$\bs y^{1;\, c_1}  = (y^{1;\, c_1}_1, \dots , y^{1;\, c_1}_r) \in \PCr$, where
\bea
y^{1;\, c_1}_{i_1}(x)\ =\ \tilde y_{i_1}(x) \ +\ c_1\, y_{i_1}^0(x)
\qquad
\text{and}\
y^{1;\, c_1}_{j}(x)\ =\ y_{j}^0(x)
\ {} \text{for}\ {}
j \neq i_1 \ .
\eea
In particular, $\bs y^{1;\,\infty} = \bs y^0$ in $\PCr$.

Assume that the family $\bs y^{d-1;\, c_1, \dots , c_{d-1}}
 \in \PCr$ is already defined.
Let \linebreak
$\tilde y_{i_d}^{\ d-1;\, c_1, \dots , c_{d-1}}$ 
be a polynomial satisfying equation
\bea
y_{i_d}^{d-1;\, c_1, \dots , c_{d-1}}(x+h) 
\ \tilde y_{i_d}^{\ d-1;\, c_1, \dots , c_{d-1}}(x)
& -&
y_{i_d}^{d-1;\, c_1, \dots , c_{d-1}}(x) 
\ \tilde y_{i_d}^{\ d-1;\, c_1, \dots , c_{d-1}}(x+h)
\phantom{ccccccccca}
\\
&=& 
 T_{i_d}(x) 
\prod_{j,\ j\neq i_d} (\,y_{j}^{d-1;\, c_1, \dots , c_{d-1}}(x+b_{i_d,j})\,)^{- a_{i_d,j}}\ .
\eea
We fix $\tilde y_{i_d}^{\ d-1;\, c_1, \dots , c_{d-1}}$
 assuming that the coefficient of $x^{l^{d-1}_{i_d}}$ in 
$\tilde y_{i_d}^{\ d-1;\, c_1, \dots , c_{d-1}}$
is equal to zero.
Set $\bs y^{d;\, c_1, \dots , c_{d}}  = (y^{d;\, c_1, \dots , c_d}_1, 
\dots , y^{d;\, c_1, \dots , c_d}_r) \in \PCr$, where
\bea
y^{d;\, c_1, \dots , c_{d}}_{i_d}(x)\ =\ \tilde y_{i_d}^{\ d-1;\, c_1, \dots , c_{d}}
(x) \ +\ c_d\,
y_{i_d}^{d-1;\, c_1, \dots , c_{d-1}}(x)
\eea
and
\bea
y^{d; \,c_1, \dots , c_{d}}_{j}(x)\ =\  
y_{j}^{d-1;\, c_1, \dots , c_{d-1}}(x)
\ {}\ \text{for} \ {}\ 
j \neq i_d \ .
\eea
In particular, $\bs y^{d;\, c_1, \dots , c_{d-1}, \infty} =
\bs y^{d-1;\, c_1, \dots , c_{d-1}}$ in $\PCr$.

The $d$-th family  is obtained from the $(d-1)$-st family by the generation procedure in
the $i_d$-th direction, see 
Section \ref{general procedure}. For any $(c_1, \dots , c_k) \in (\Bbb P^1)^k$ the tuple
$\bs y^{k;\, c_1, \dots , c_{k}}$ lies in $P$.

For any $(c_1, \dots , c_k) \in \C^k$ and any $i \in \{1, \dots , r\}$, we have
\bea
\deg\
y^{k; \,c_1, \dots , c_{k}}_{i}(x)\ =\  l^w_i\ .
\eea
Set
\bea
P^{(i_1, \dots , i_k)}\ = \
\{ \ \bs y^{k;\, c_1, \dots , c_{k}} \ |\
(c_1, \dots , c_k) \in \C^k\ \} \ .
\eea

For every $(c_1, \dots , c_k) \in (\bP)^k$ we define a rational function
$v_{c_1,\dots , c_k} : \C \to N_+$ by the formula
\bean\label{S}
v_{c_1,\dots , c_k} : x \mapsto
\exp\,(g_k(x; c_1, \dots , c_k) E_{i_k})
\dots
(\exp\,(g_1(x; c_1) E_{i_1})
\eean
where
%\bea
%g_1(x; c_1)
% =  \frac{ T_{i_1}( x ) \prod_{j,\ j\neq i_1} (\,y^0_j(x+b_{i_1,j})\,)^{- a_{i_1,j}}}
%{y_{i_1}^{1;\,c_1}(x)\ y^0_{i_1}(x)} ,
%\eea
\bea
g_d(x; c_1, \dots , c_d) = 
\frac{ T_{i_d}(x)  
\prod_{j,\ j\neq i_d } (\, y_{j}^{d-1;\, c_1, \dots , c_{d-1}}(x+b_{i_d,j}) \,)^{- a_{i_d,j}}}
{ y_{i_d}^{d;\, c_1, \dots , c_d}(x)\  y_{i_d}^{d-1;\, c_1, \dots , c_{d-1}}(x)}
\eea
for $d = 1, \dots , k$. 
In particular, if some of $c_1, \dots, c_k$ are equal to $\infty$, then 
the corresponding exponential factors in \Ref{S} must be replaced by the 
identity element $\id \in G$.

The function $v_{c_1,\dots , c_k}$ continuously depends on 
$(c_1, \dots , c_k) \in (\Bbb P^1)^k$. For any $(c_1, \dots , c_k) \in (\Bbb P^1)^k$
the pair $(D^{v_{c_1,\dots , c_k}}, v_{c_1,\dots , c_k})$ belongs to
$\widehat\tN$.

Let $x_0\in \C$ be a regular point. Consider the map
\bea
 \phi :         \C^k \to G/B_-,
\qquad
(c_1, \dots , c_k) \mapsto (v_{c_1,\dots , c_k}(x_0))^{-1}B_- \ .
\eea

\begin{proposition}\label{MAIN}
The image of the map $\pi$  is $B_{w^{-1}}$.
\end{proposition}

\noindent
{\it Proof of Proposition \ref{MAIN}.}
For any $(c_1, \dots , c_k) \in (\bP)^k$ consider the rational section
\bean\label{SECT}
 S_{(c_1, \dots , c_k)} \ :\ x \ \mapsto\ 
( (v_{c_1,\dots , c_k}(x))^{-1} B_-)\times x 
\eean
of the bundle $p'$.
In particular, if some of $c_1, \dots, c_k$ are equal to $\infty$, then 
the corresponding exponential factors in \Ref{S} must be replaced by
the identity element  $\id \in G$.
This section is horizontal with respect to the connection
$\nabla_{\bs y^0}'$ and continuously depends on 
$(c_1, \dots , c_k) \in (\Bbb P^1)^k$.
This means that
\bea
\phi (\C^k) \ \subset \ B_{w^{-1}}\ ,
\eea
see Corollary \ref{Bruhat},
and if some of $c_1, \dots , c_k$ are equal to $\infty$, then
$S_{(c_1, \dots , c_k)}(x) \notin B_{w^{-1}}$, see Corollary \ref{Bruhat-2}.
It remains to show that every point in $B_{w^{-1}}$ is the limit of points
of the form $S_{(c_1, \dots , c_k)}(x)$. But that statement follows from
Corollary \ref{Bruhat} and

\begin{lemma} 
Assume that $x_0\in \C$ is such that $T_i(x_0)\neq 0$, $y^0_i(x_0)\neq 0$, for
$i = 1, \dots , r$. Assume that $x_0\in \C$ is such that 
$y^0_j(x_0 + b_{i,j}) \neq 0$ for all $i\neq j$.
Then there exists a proper algebraic subset $K \subset (\C - 0)^k$ with 
the following property. For any 
$(c_1^1, \dots , c_k^1)\, \in \, (\C - 0)^k - K$ there exists
a unique $(c_1^2, \dots , c_k^2) \in \C^k$ such that
\bea
(c_1^1, \dots , c_k^1)\ =\ (g_1 (x_0; c_1^2), \dots , 
g_k ( x_0; c_1^2, \dots , c_k^2)) \ .
\eea
\hfill
$\square$
\end{lemma}

The proposition is proved.
\hfill
$\square$

Theorem \ref{second main} is a direct corollary of Proposition \ref{MAIN}.

\end{document}